\documentclass[12pt]{article}
\usepackage{mathpazo}
\usepackage{stmaryrd}
\usepackage[leqno]{amsmath}
\usepackage{amssymb,theorem}
\begin{document}
\title{Generic structures\thanks{I owe much to the members of the FSB seminar in Bristol, the audience in the Philosophy of Mathematics seminar in Oxford (especially James Studd), Hazel Brickhill, Kit Fine, Philip Welch, \O ystein Linnebo, and Giorgio Venturi. I am deeply grateful to the organisers (David Svoboda and Prokop Sousedik) of the conference on \emph{The emergence of structuralism and formalism} (Prague, 2016), where I gave my first talk on this subject. I am also very grateful to the members of the audience in Prague for helpful questions and suggestions. To conclude, I am indebted to two anonymous referees for their comments on an earlier version of this article.} }  
\author{Leon Horsten \\ Department of Philosophy \\ University of Bristol \\ Bristol, BS6 6JL \\ United Kingdom \\ Leon.Horsten@bristol.ac.uk \\  https://orcid.org/0000-0003-3610-9318
}

\date{}

 \maketitle

\begin{abstract}
\noindent In this article ideas from Kit Fine's theory of arbitrary objects are applied to questions regarding mathematical structuralism. I discuss how sui generis mathematical structures can be viewed as generic systems of mathematical objects, where mathematical objects are conceived of as arbitrary objects in Fine's sense.
 \end{abstract}

%%%%%%%%%

\section{What are mathematical theories about?}

Many philosophers today consider mathematical \emph{structures} to be the subject matter of mathematics. 
On the one hand there is \emph{sui generis} or \emph{non-eliminative} structuralism. According to sui generis structuralism, the subject matter of a mathematical theory is a \emph{mathematical structure} or a family of mathematical structures, where mathematical structures are understood to be \emph{abstract universals}. On the other hand there is \emph{eliminative} structuralism, which is not wedded to abstract universals. %Eliminative structuralist accounts hold that when I informally say that two systems of objects instantiate the same structure, all this means is that the systems in question are isomorphic.

Inspired by work of Fine, I explore in this article a new answer to the question what mathematical structures are. My account makes use of elements of Fine's theory of arbitrary or generic objects \cite{Fine 1983}, \cite{Fine 1985}, and extends this to a theory of arbitrary or generic \emph{systems} of objects \cite{Fine 1998}. Thus it combines elements of object-Platonism with elements of structure-Platonism. 

Following a tradition in the discussion of structuralism in the philosophy of mathematics, I  take the distinction between \emph{algebraic} mathematical theories  (such as group theory) and \emph{non-algebraic} mathematical theories (such as arithmetic)  to be important in this context. On the proposed account, 
a non-algebraic mathematical theory is about \emph{arbitrary} (or generic) \emph{objects}. These arbitrary objects are abstract, and they form an abstract \emph{structure}, which is a generic system of arbitrary objects. An algebraic theory, in contrast, describes not one single structure but about a \emph{family} of structures, where again a structure is understood as a generic system of objects. 
%So generic structuralism is a Platonistic theory of both algebraic and non-algebraic theories.

I proceed as follows. First, I state my reasons for being dissatisfied with existing versions of both non-eliminative and eliminative structuralism (section \ref{mathstruc}). Next, I discuss Kit Fine's theory of arbitrary objects and how it can be extended to a theory of generic systems (section \ref{arbitrary objects}, section \ref{gen}) and to an interpretation of the theory of pure sets (section \ref{sets}). This will lead to a new account of the nature of mathematical structures; I will explain in some philosophical detail what this position amounts to  (section \ref{gen-struc}). Before closing, I compare my position with rival accounts (section \ref{comp}).

%I WILL TRY TO ARTICULATE A NEW POSITION; IT IS GOING TO BE SKETCHY

%%%%%%%%%%%%%

\section{Mathematical structuralism}\label{mathstruc}

%Plato famously took Science, including Mathematics, to be about \emph{Ideas} (or \emph{Forms}). He took these Ideas to be abstract and to exist independently of the concrete world that we inhabit. 
%Aristotle disagreed with Plato about the independent existence of Ideas: he maintained that Forms only ever exist embodied in concrete objects. 
%Aristotle held that our concrete world is in some fundamental sense always finite. This presents a problem for the philosophy of mathematics, for mathematics is concerned with Ideas that are in some sense infinite (such as the idea of Number). The concept of potential infinity was developed to accommodate mathematics in an ever finite world.
%In the twentieth century, the concept of Idea or Form was replaced by the concept of \emph{structure}. Thus the revised thesis becomes that Science, including Mathematics, is about structures. 

A distinction is drawn between eliminative structuralism and non-eliminative structuralism \cite[p.~81]{Shapiro 1996}. 
Non-eliminative and eliminative structuralists agree that different \emph{systems} of concrete objects can have a structure in common. But Platonistic structuralism maintains that such mathematical structures exist independently of the systems that \emph{instantiate} them; eliminative structuralists regard talk of mathematical structures as loose talk that can ultimately be replaced by talk about systems being isomorphic to each other.

Another distinction can be drawn between \emph{algebraic} and \emph{non-algebraic} mathematical theories \cite[p.~40--41]{Shapiro 1997}. Intuitively, an algebraic theory is one that is intended to be about many different structures. A non-algebraic theory, in contrast, intends to describe one structure only.

Eliminative structuralists hold that \emph{every} mathematical theory, algebraic or non-algebraic, is about a multiplicity of systems of objects. The non-eliminative structuralist can only partly agree with this thesis of distributed reference for mathematical theories. She agrees with the eliminative structuralist that algebraic theories are about many structures. But she insists that the distributed reference hypothesis does not hold for non-algebraic theories and claims to have evidence for this: we speak of ``\emph{the} natural number structure''. Non-algebraic theories are about a unique subject matter: the subject matter of a non-algebraic theory is the unique \emph{structure} that it purports to describe.\footnote{In more recent work, Shapiro questions whether non-eliminative structuralism should commit itself to this unique reference thesis: see \cite[p.~243]{Shapiro 2006}.}

Eliminative structuralism comes in many flavours. One popular variant is \emph{set theoretic structuralism}, which takes the structures that a mathematical theory is about to be \emph{sets} endowed with operations \cite{Mayberry 1994}. But this assumes a form of set-theoretic reductionism that many today find hard to accept. 
%Even though groups, for instance, can be taken to be sets ---indeed, the first pages of a typical textbook on group theory will define it as such,--- this fact will normally be of little or no relevance to a group theorist, who treats groups as if they exist independently of the universe of sets.
Eliminative structuralism does not \emph{have} to take structures to be sets. She can, for instance, take structures to be physical arrangements (``pluralities'') of objects that stand in specific physical relations to each other. But this is taking mathematical structures to be physical in nature, which is is another form of reductionism that is hard to accept. I will not pursue this debate here.

More importantly for present purposes, according to eliminative structuralism there is no sense in which even non-algebraic theories can be said to be about mathematical \emph{objects}. But mathematical theories are about mathematical objects \cite[p.~1]{Parsons 2008}:
\begin{quote}
The language of mathematics speaks of objects. This is a rather trivial statement; it is not clear that we can conceive any developed language that does not. What is of interest is that, taken at face value, mathematical language speaks of objects distinctively mathematical in character: numbers, functions, sets, geometric figures, and the like. To begin with, they are distinctive in being abstract.
\end{quote}

Existing forms of non-eliminative structuralism cannot be charged with implausible reductionist claims. Moreover, according to non-eliminative forms of mathematical structuralism, non-algebraic theories can be taken to be about objects. Consider Shapiro's version of non-eliminative structuralism \cite{Shapiro 1997}. On this view, a mathematical structure contains places or roles that can be \emph{occupied} by objects. Nonetheless, the places themselves can be viewed as objects that can be organised into a system that instantiates the structure   \cite[p.~100--101]{Shapiro 1997}.

However, at this point the non-eliminative structuralist is faced with two problems.

First, the objects that populate \emph{sui generis} structures are according to non-eliminative structuralism in some sense \emph{incomplete}. A rough statement of what is intended (when applied to arithmetic) is to say that numbers only have \emph{structural properties}  \cite[p.~72--73]{Shapiro 1997}. However, it has turned out to be very difficult to make the intended meaning of such statement sufficiently precise in a way that does not lead to counterexamples.\footnote{See \cite{Linnebo and Pettigrew 2014} and \cite{Korbmacher and Schiemer forthc.}.}
 Shapiro has also come to recognise that it is difficult to give a satisfactory philosophical account of the incompleteness of mathematical objects \cite[section 1]{Shapiro 2006}. For instance, the number 7 might have the property of being my least favourite number, even though this is not a `structural' property.

Secondly, as Hellman points out \cite[p.~546]{Hellman 2006}, Shapiro's view is vulnerable to a \emph{permutation objection}. If, in the \emph{ante rem} structure $\mathbf{N}$ of the natural numbers, we permute its places (in a non-trivial way), then we obtain a system $\mathbf{N'}$ that is isomorphic to $\mathbf{N}$. We can then ask, in the spirit of \cite{Benacerraf 1965}: what could possibly make $\mathbf{N}$, rather than $\mathbf{N'}$ be the unique \emph{sui generis} structure that arithmetic is about?

The aim of this article is to articulate and explore a new position according to which mathematical theories are \emph{not only about structures but also about objects}. So it cannot be a form of eliminative structuralism. The position to be developed is expected to attribute  the `right' kind of incompleteness to natural numbers, and should not be vulnerable to the permutation objection.

\section{Arbitrary objects}\label{arbitrary objects}

A theory of arbitrary objects was developed in \cite{Fine 1983} and \cite{Fine 1985}. I will now briefly and informally review some the main tenets of this theory.\footnote{In what follows, I may be guilty of adding a few tenets of my own to Fine's theory.}

Consider any category  of entities. It is helpful to fix (without essential loss of generality) on some particular kind of mathematical objects: the natural numbers, say. There are \emph{specific} natural numbers, such as the number 23. But beside specific natural numbers, there are also \emph{arbitrary} natural numbers. An arbitrary number is what a mathematician refers to when she says: ``Let $m$ be a natural number\ldots'',\footnote{This is denied by \cite{Breckenridge and Magidor 2012}: see section \ref{magidor} below.} and then goes on to reason about $m$. 

There are many arbitrary natural numbers. For instance, it would make perfect sense for our mathematician, in the course of her argument, to add ``Now let $n$ be an \emph{other} natural number\ldots'', and go on reasoning about both $m$ and $n$.

In general, an arbitrary natural number does not determinately have any specific natural number as its value. There is no determinate matter of fact, for instance, about whether the value of our mathematician's arbitrary number $m$ is 23.\footnote{Again, \cite{Breckenridge and Magidor 2012} deny this.}

There \emph{can} be a determinate fact about whether an arbitrary number $x$ is numerically identical with an arbitrary number $y$. Our mathematician was perfectly within her rights when she required the arbitrary numbers $m$ and $n$ to be non-identical. She might just have said, more clearly perhaps: ``take any two arbitrary numbers $m$ and $n$ such that $m \neq n, \ldots$''

When an arbitrary natural number does not determinately have some given specific number as its value, there is a sense in which it \emph{can} be the specific number in question. Thus we can say that arbitrary numbers can be in different specific \emph{states}. These possible situations (states) may be understood in a Lewisian realist way, but one can also understand talk about states in a Kripkean, deflationist way. (I will not take a stance on this matter.)

There is, however, no \emph{actual} specific state in which the arbitrary number is. The best we can say, perhaps, is that an arbitrary number ``actually'' is in a ``superposition'' of specific states.

There are \emph{degrees} of arbitrariness. If our mathematician says ``Let $m$ be an arbitrary natural number larger than 10'', then she refers to a number that is less arbitrary than when she says ``Let $m$ be an arbitrary natural number larger than 5''. So one might say that a natural number is \emph{completely arbitrary} if it can be any specific natural number whatsoever. In some sense, specific natural numbers are closely related to \emph{limiting cases} of arbitrary natural numbers. For many purposes they can be identified with arbitrary natural numbers that have an absolutely minimal degree of arbitrariness. But it is important to see that, strictly speaking, \emph{no} arbitrary number is identical to any specific number. An arbitrary number is the sort of thing that can be in a state, whereas a specific number cannot be in a state.

Since a completely arbitrary natural number can be in any specific state whatsoever, a completely arbitrary natural number has only those `specific' (Fine calls them `classical') properties that every specific natural number has. Let us call this \emph{Fine's principle}.  It is not a simple matter to spell out precisely what `specific' properties are,\footnote{Fine recognises that it is difficult to give a precise description of the distinction between `classical' and `generic' conditions \cite[chapter 1]{Fine 1985}. For a discussion of this distinction, see \cite[section 2.1.3]{Breckenridge and Magidor 2012}. } and I will not attempt to do so here. But self-identity is one such specific property of the completely arbitrary number $m$ above, whilst ``possibly being identical with the arbitrary number $n$ (above)'' isn't.

An arbitrary number can be more or less \emph{likely} to be in a given state. For instance, an arbitrary natural number between 25 and 50 is unlikely to be a power of 2.\footnote{Fine does not connect, as I do, the concept of arbitrary object with a concept of probability. Since the relevant probability functions should be uniform distributions on an infinite space, it might be appropriate to appeal to techniques of \cite{Benci et al 2013} to model them. However, I leave this discussion for another occasion.}

The concept of arbitrary object has only played a marginal role in philosophy. There are two reasons for this. First, arguments have been put forward that purport to show that it is philosophically untenable to maintain that there are arbitrary objects \cite[p.~160]{Frege 1979}. Second, it has been argued that the concept of arbitrary object cannot be put to \emph{good use}. For instance, Russell's doctrine of incomplete symbols shows how it is not necessary to take the phrase ``the man in the street'' to be a denoting term referring to some `arbitrary man'. 
Fine intended to counter both of these arguments. 
%It does not suffice to counter only the first: even if a natural and coherent theory of arbitrary objects can be worked out, this effort is no more than an academic exercise if there is no essential use for them. And this is where, perhaps, clear success has so far been lacking.

Fine sought to defuse the first objection by articulating a \emph{coherent} and \emph{natural} concept of arbitrary object. I assume for the purposes of the discussion that he was successful in this enterprise. 

Fine countered the second objection in a double movement. On the one hand, he gestured at a number of other applications for the concept of arbitrary object, such as the theory of infinitesimals, the notion of forcing in set theory,\ldots\footnote{As far as I know, these potential applications have not yet been worked out in detail.} On the other hand, he worked out two applications in detail. First, he formulated a natural semantics for first-order logic in terms of the notion of arbitrary object \cite{Fine 1985}. Second, he discussed in detail how a theory of Cantorian (and Dedekindean) abstraction could be developed based on the theory of arbitrary objects \cite{Fine 1998}.

Russell showed that many uses of arbitrary objects are \emph{non-essential}: much talk involving arbitrary objects can be adequately paraphrased in ways that do not involve them. A somewhat similar point can be made about Fine's use of arbitrary objects in his new semantics for first-order logic. The silence with which his effort seems to have been met\footnote{An exception to this is \cite{King 1991}.} is due to the fact that we have perfectly adequate semantics for first-order logic that do not involve arbitrary objects.\footnote{The situation is somewhat similar, in this respect, to theories of the mathematical continuum that involve infinitesimals.}

Fine's account of Cantorian abstraction, however, does strike me as a use of the theory of arbitrary objects that is much harder to dismiss as philosophically redundant. Moreover, Fine himself observed that this account can be extended to form the basis of a new form of mathematical structuralism for non-algebraic mathematical theories. A primary aim of the present article is to explore the ramifications of Fine's brief but suggestive remarks in section VI of \cite{Fine 1998} on these matters in some depth, not only for the interpretation of non-algebraic but also for the interpretation of algebraic theories. My account will diverge in key points from the way in which Fine thought that this new form of structuralism should look like. I will point out where this is the case, and I will state my reasons for developing it in a different way.

%%%%%%%%

\section{Generic systems}\label{gen}

Generic systems can be seen as a special kind of arbitrary entities. My account of generic systems is intended to be completely general. In order to illustrate how it works for particular mathematical theories, I will concentrate on one non-algebraic theory (arithmetic), and one algebraic theory (graph theory). 

I take what I am doing to be an exercise in \emph{naive metaphysics} in the sense of \cite{Fine 2017}. That is, I investigate the \emph{nature} of generic systems. As a tool in this investigation, I explore in this section how some generic systems can best be modelled in set theory. The purpose of this is to discern metaphysical properties of generic systems. For instance, a set theoretic way of modelling a generic system will give an answer to questions such as: How is a generic system incomplete? How many objects does it contain? What kinds of things are the states? In how many states can this generic system be? Such answers are only as good as the set theoretic model is at representing fundamental properties of the generic system in question. But generic systems, like Fine's arbitrary objects, constitute a metaphysical realm in their own right: far be it from me to advocate an ontological reduction of generic systems to sets.

%My account is to some extent inspired by \cite{Urquhart 2015}.

%%%%%%%%%%

\subsection{Maximally arbitrary $\omega$-sequences}\label{arbitrary omega sequences}

Suppose that the physical world around us consists of a countably infinite collection $A$ of objects, which we may label as $a_0, a_1,a_2 \ldots$ 

Consider all ways in which this collection can be ordered as an $\omega$-sequence. Mathematically, this amounts to considering all functions from $\mathbb{N}$ to $A$. Let us write down all such functions in a long list $\mathrm{L}$. It is clear that this list has $2^{\aleph_0}$ entries: let $p_{\alpha}$ (for every $\alpha < 2^{\aleph_0}$) be line $\alpha$ in this long list $\mathrm{L}$.

The entries in this long list (i.e., the $p_{\alpha}$'s) constitute all the possible \emph{specific} $\omega$-sequences in the physical world that we are considering. Thus the $p_{\alpha}$'s can be seen as exhausting the possible \emph{states} that a generic $\omega$-sequence can be in.

The possible state space should be sufficiently large to accommodate every possible specific state. But there seems no reason to assume that it is larger. So we propose that the $p_{\alpha}$'s are labels for the possible states. Then we can take the list $\mathrm{L}$ to give the \emph{possible state profile} of an arbitrary $\omega$-sequence. The state profile of a generic  $\omega$-sequence is all we care about, so we may for modelling purposes take the list $\mathrm{L}$ to \emph{be} an arbitrary $\omega$-sequence. 

Of course there are infinitely many other such lists of order type $2^{\aleph_0}$ ---there are in fact $2^{2^{\aleph_0}}$ of them. We take each of them to specify the \emph{state profile} of some arbitrary $\omega$-sequence, or, in brief, to \emph{be} an arbitrary $\omega$-sequence.

The list $\mathrm{L}$ is a \emph{maximally arbitrary} $\omega$-sequence. A much less arbitrary $\omega$-sequence is given by a list (of cardinality $2^{\aleph_0}$) of functions from $\mathbb{N}$ to $A$ almost all of which are identical to some particular function $f$. This arbitrary $\omega$-sequence is then overwhelmingly likely to be the specific $\omega$-sequence $f$. And specific $\omega$-sequences can then be seen as canonically embedded in the arbitrary $\omega$-sequences. They are embedded as the lists that have the same $\omega$-sequence on each row.

Generic systems consist of arbitrary objects. Generic $\omega$-sequences consist of arbitrary natural numbers. To see how this works, consider again our long list $\mathrm{L}$. An arbitrary natural number is a \emph{thread} or \emph{fiber} through $\mathrm{L}$, i.e., formally, a function from $2^{\aleph_0}$ to $A$. So there are $2^{2^{\aleph_0}}$ arbitrary natural numbers in $\mathrm{L}$. Elementary arithmetical operations on arbitrary natural numbers are defined pointwise. Take the sum $a+b$ of the arbitrary numbers $a$ and $b$, for instance. This is the arbitrary number that in each state $w$ takes the value of the sum of the value of $a$ at $w$ and the value of $b$ at $w$ according to the state ($\omega$-sequence) $w$. Elementary operations can then straightforwardly be seen to satisfy the familiar properties of arithmetical operations (such as commutativity of $+$, for instance).\footnote{The mathematical properties of the generic $\omega$-sequence are investigated in \cite{Horsten and Speranski 2018}.}

%Now consider, for every $n \in \mathbb{N}$, the function $f_n$ which takes as its function value, for every $\alpha < 2^{\aleph_0}$, the element of $A$ that appears as the $n$-th entry in $p_{\alpha}$. Each such $f_n$ is a thread through $\mathrm{L}$, i.e.,  an arbitrary natural number. Then the list $\mathrm{L}$  can be seen as an $\omega$-sequence of arbitrary numbers: it can be parsed as the $\omega$-sequence $f_0, f_1, f_2, \ldots$

%It follows from this that when I informally said in section \ref{arbitrary objects} that arbitrary numbers \emph{can be} specific numbers, I was speaking ambiguously. If a specific number is a limiting case of an arbitrary number, then it cannot literally \emph{be} a (non-minimally) arbitrary number; but it can take on the same value as a non-minimally arbitrary number.

%%%%%%%%

\subsection{The generic countable graph}\label{arbitrary countable graph}

Let us do this also for countable simple graphs, and let us do it in the same way. Informally it is clear what we should do, but, for definiteness, let us see in detail how it goes.

Suppose we are given a countable vertex set $V=  v_0,v_1,v_2,\ldots$ We want to list all possible \emph{particular} graphs on subsets of $V$. 

A graph consists of vertices and edges. Edges are often seen as unordered pairs of vertices. But in the framework of multisets, vertices themselves can be seen as special cases of edges. So let us take this perspective.
Let an edge between two vertices $v_i$ and $v_j$ be given as the unordered (multi-)pair $\{ v_i,v_j \}$. And let a vertex $v_i$ be given as the multipair $\{ v_i,v_i \}$.

We order the unordered multi-pairs of elements of $V$ in lexicographical fashion, i.e., when $i \leq j$ and $i'\leq j'$, we say that \begin{quote}$\{ v_i,v_j \} < \{ v_{i'},v_{j'} \}$ if and only if $i<i'$ or ($i=i'$ and $j  < j'$).\end{quote} Then a graph $G$ is given by a countable list of unordered multi-pairs, without repetitions, well-ordered by $<$, where we adopt the convention that if $\{ v_i,v_j \}$ appears in the list, the multipairs $\{ v_i,v_i \}$ and  $\{ v_j,v_j \}$ must also belong to the list.

Consider a list $\mathrm{L}^*$ of such graphs which contains every particular graph on $V$ exactly once. This list will again be of length $2^{\aleph_0}$. The list $\mathrm{L}^*$ gives the state profile of a \emph{maximally arbitrary countable graph} on $V$, where a row in $\mathrm{L}^*$ describes a specific state that this graph can be in.
So we can take it, for modelling purposes, to \emph{be} a maximally arbitrary countable graph. There are $2^{2^{\aleph_0}}$ such lists, i.e., generic countable graphs on $V$. 

Generic graphs then consist of arbitrary edges. An arbitrary edge of the arbitrary graph $\mathrm{L}^*$ is a thread through $\mathrm{L}^*$. So there are $2^{2^{\aleph_0}}$ arbitrary edges in $\mathrm{L}^*$. 
%It is then easy to see how $\mathrm{L}^*$ can be seen as a list of arbitrary edges.

Some generic graphs are less arbitrary than others. There are lists of graphs (of length $2^{\aleph_0}$) that list the same particular graph on almost all of their rows. Such arbitrary graphs are overwhelmingly likely to be one particular graph. And there are generic graphs that list the same specific graph on each of its rows, so particular graphs are canonically embedded as limiting cases of generic graphs. 

A celebrated theorem from countable graph theory (by Erd\"{o}s \& Renyi, 1963) says that there exists a simple graph $R$ with the following property \cite{Cameron 2013}. If a countable
graph is chosen at random, by selecting edges independently with probability $\frac{1}{2}$ from the set of two-element subsets of the vertex set, then almost surely (i.e., with
probability 1), the resulting graph is isomorphic to $R$. This graph $R$ is called the \emph{Rado graph}. In the present context, this theorem can be taken to say that with probability 1, the generic countable graph $L^*$ is the Rado graph. Thus Erd\"{o}s \& Renyi's theorem can be seen as a \emph{probabilistic categoricity theorem} for countable graph theory.

We could now do the same for some other algebraic theory, such as countable group theory for instance. We would then take a particular countable group to be given by its multiplication table, and take a generic countable group to be given by a long list of particular countable groups. We could also move from countable mathematics to uncountable mathematics, by considering the theory of the real numbers for instance. But the general procedure is clear, so there is no need to go into details here. 

%%%%%%%%%%

\section{Sets and arbitrary objects}\label{set theory}\label{sets}

Nothing has yet been said about the nature of the objects of the underlying domain of generic systems ($\omega$-sequences, graphs, groups): the choice of the underlying plurality of objects in section \ref{gen} was an external parameter of the model. But this is a question that cannot altogether be avoided.

Nothing in mathematical practice dictates that the underlying objects should be physical in nature. But clearly \emph{many} of them are required if we want to extend the view described in this article beyond countable mathematics.

The \emph{universe of sets} contains a sufficient supply of objects for our mathematical objects and structures.  It might seem \emph{ad hoc} if the subject matter of set theory would fall altogether outside the scope of the theory of mathematical structures.

If we admit higher order arbitrary objects, then the theory of arbitrary objects provides the resources for giving a natural explication of what an iterative hierarchy of set theoretic ranks may be taken to consist of. A metaphysical account of how this might work goes along the following lines:

\begin{itemize}

\item \textbf{Stage 0}

We take some specific object $o$ to be given.

\item \textbf{Stage 1}

Now consider the arbitrary object which can only be in the state of being the specific object $o$; we denote this arbitrary object as $\langle o \rangle$. Then $o \neq \langle o \rangle$: unlike an arbitrary object, a specific object is not the sort of entity that can be in states.

\item \textbf{Stage 2}

(2a) Consider the higher-order arbitrary object that can only be in the state of being the arbitrary object $\langle o \rangle$. Denote this higher-order arbitrary object as $\langle \langle o \rangle \rangle$. 

(2b) Consider an higher-order arbitrary object that can only be in one of the following two states: being the specific object $o$, or being the arbitrary object $\langle o \rangle$. Denote this object as $\langle o, \langle o \rangle  \rangle$.

\item \textbf{Stage $\omega$}

Collect the arbitrary objects that have been generated in the finite stages.

\item \textbf{Later stages}

\emph{Continue in this way into the transfinite.}

\end{itemize}

\noindent In this way, a hierarchy of higher-order arbitrary objects, based on one specific object, is built up.
The simple idea of course is to view arbitrary objects as (non-empty) sets, and their states as their elements. This leaves the specific object on which the hierarchy is based, which is not the sort of thing that can be  in states, playing the role of the empty set. So we read the pointy brackets $\langle \ldots \rangle$ as curly brackets $\{ \ldots \}$, and $o$ as $\emptyset$. For example, the set $\{ \emptyset, \{ \emptyset \} \}$ is the arbitrary object which can be specific object $o$, but which can also be the arbitrary object which can only be the specific object $o$.

Let us call this the \emph{generic hierarchy}. The concept that it might \emph{perhaps} be taken to capture is the \emph{combinatorial} set concept, which takes a set of $X$-es to be the result of an arbitrary selection process.

The construction of the generic hierarchy makes use of a notion of higher order arbitrariness, which in turn requires taking arbitrary objects ontologically very seriously. So this is a decidedly Platonistic interpretation of set theory.

We have seen how the distinction between $a$ and $\langle a \rangle$ corresponds to the distinction between a set and its singleton. It is not clear whether this metaphysical account gives an adequate explanation of what Lewis has dubbed the ``mystery of the singleton relation'' \cite{Lewis 1991}. Some may worry that the distinction between a specific object and the arbitrary object that can only ever be in the state of being that specific object is no better motivated than the relation between an object and its singleton.

An important observation is that the generic hierarchy is the result of an extreme \emph{flattening} or extensionalisation of what is at bottom a much more intensional notion of set. There are, for example, in fact \emph{two} higher order arbitrary objects that can play the role of the set $\{ \emptyset, \{ \emptyset \} \}$ in stage (2b). These two arbitrary objects are anti-correlated; if one of them is in the state of being $o$, then the other is in the state of being $\langle o \rangle$ and vice versa. (So one might denote these two objects as $\langle o, \langle o \rangle  \rangle$ and $\langle \langle o \rangle, o  \rangle$, respectively.) 
%Moreover, in analogy with before, it is not hard to see how in a very natural way arbitrary sets arise which are \emph{more likely} to be in one state than in another.

The engine of the iterative hierarchy of sets is the full power set operation. This operation drops out naturally in the framework of the theory of arbitrary objects: the description of stages shows how at every stage $\alpha + 1$, there are $2^{\alpha}$ arbitrary objects that naturally represent all the subsets of the domain of sets that exist at stage $\alpha$.\footnote{I do not claim that the present account has much new to say about the motivation of some of the other powerful axioms of standard set theory such as Infinity or Replacement.}

%%%%%%

\section{What mathematical theories are about}\label{gen-struc}

Let us now connect the discussion of generic systems to the question of the subject matter of mathematical theories.

%%%%%%%%%%%%

\subsection{Non-algebraic theories}\label{non-alg}

In the setting of arbitrary objects and systems, sense can be made of the commonplace statement that arithmetic is about the natural numbers. In fact, there are two \emph{different} ways in which this can be made more precise. One might say that arithmetic is about all \emph{specific} natural numbers, where specific numbers are arbitrary objects as described in section \ref{arbitrary omega sequences}. Alternatively, one might say that arithmetic is about \emph{arbitrary} natural numbers.\footnote{The form of structuralism sketched in \cite[section 6]{Fine 1998} takes arithmetic to be about the \emph{specific} natural numbers rather than about the specific and the arbitrary natural numbers.} But in either case, arithmetic is said to be about \emph{objects} of an abstract kind, in accordance with object-focussed forms of mathematical Platonism such as G\"{o}del's. I will not try to decide here whether it is more natural to say that arithmetic is about arbitrary natural numbers or about specific natural numbers.

As an illustration, consider \emph{Goldbach's conjecture}, which says: 
\begin{quote} Take any even natural number larger than 2; \\ it is the sum of two prime numbers. \end{quote} Two possible readings are:
\begin{enumerate}
\item Take any \emph{specific} even natural number $n$; then $n$ is the sum of two (specific) primes.
\item Take any \emph{arbitrary} even natural number $n$; then $n$ is the sum of two (arbitrary) primes.
\end{enumerate}

From the point of view of our set theoretic model, the arbitrary natural numbers are located (as threads) in a large structure: the collection of all $\omega$-sequence orderings of our countable collection $A$. So there is a sense in which arithmetic can be said to be about the generic $\omega$-sequence. 

The generic $\omega$-sequence can be interpreted in a \emph{sui generis} manner. 
Alternatively, one can `deflate' the generic $\omega$-sequence in an eliminative way as supervening on concrete $\omega$-sequences. But then one cannot take arbitrary objects (i.e., specific and arbitrary numbers) to exist as completed abstract entities, and that would mean that arithmetic is not about mathematical \emph{objects}. So generic $\omega$-sequences are to be understood in a non-eliminative way.

We have seen in section \ref{arbitrary omega sequences} that there are, from a higher-order point of view, \emph{many} (completely) arbitrary $\omega$-sequences, which can be obtained from each other by permutations of the list $\mathrm{L}$. 
%So it might seem that Hellman's permutation objection still looms.
But from within arithmetic as an independent and self-standing discipline there is only \emph{one} generic $\omega$-sequence. All attempts to `construct' other ones are, from the point of view of arithmetic itself, mere re-labelings of the one and only completely arbitrary $\omega$-sequence. So from \emph{this} point of view, our mathematical modelling of the generic $\omega$-sequences as a list  contains excess structure. Therefore a better way of modelling is obtained when an arbitrary number is taken to be the set
$$ \{ \langle  w, a \rangle \mid w \textrm{ is an } \omega\textrm{-sequence on } A, \textrm{ and } a \in A     \},  $$ where $A$ is again the countably infinite set of objects of section \ref{arbitrary omega sequences}.\footnote{This is roughly how Fine understands arbitrary natural numbers (see \cite[p.~630]{Fine 1998}), except that I am here giving a set theoretic representation of arbitrary natural numbers and the structure to which they belong, whereas Fine does not do so. He instead appeals to an intuitively given theory of arbitrary objects and locates the generic structures and objects in the ontology of such a theory.}
%These arbitrary numbers together form \emph{the} generic $\omega$-sequence.

Whether this way of modelling the generic $\omega$-sequence brings out most of its structural features depends on wider issues in the philosophy of arithmetic. Some hold that we must be able to \emph{compute} on the particular $\omega$-sequences (states) that instantiate the natural number structure \cite{Halbach and Horsten 2005}, \cite{Horsten 2012}. If that is so, then the states are `computable' $\omega$-sequences, i.e., \emph{recursive} permutations of a recursive notation system for the natural numbers (our system of arabic numerals, for instance,  or the stroke notation system for the natural numbers). On this picture, the generic natural number structure can only be in $\aleph_0$ many states, and if we go on to define arbitrary numbers as before, then there will only be $2^{\aleph_0}$ many of those. Again, I will not try to adjudicate here which way of modelling the generic natural number structure is most faithful to its metaphysical nature.

At any rate, the specific natural numbers form a sub-structure of this larger structure. 
The specific natural numbers, ordered in the natural way, of course \emph{themselves} form an $\omega$-sequence. \emph{This} $\omega$-sequence is what Fine regards as ``the generic $\omega$-sequence \cite[p.~618--619]{Fine 1998}.\footnote{Fine does not discuss the further phenomenon that from a higher-order perspective, there are \emph{many} arbitrary $\omega$-sequences.} This means that on his view, in contrast to mine, the generic $\omega$-sequence \emph{instantiates itself} \cite[p.~621]{Fine 1998}. But that of course leaves the resulting account open to the permutation objection. Fine himself is well aware of this \cite[p.~625]{Fine 1998}, and he defends his account against this charge \cite[p.~628--629]{Fine 1998}. It would take me too far to go into details here, but I find Fine's arguments in this matter unconvincing. Fine's reply to the permutation objection involves, in the final analysis, making somewhat arbitrary choices and consequently introducing excess structure.

%%%%%%%%

\subsection{Algebraic theories}

We could say that countable graph theory is about the generic countable graph as modelled in section \ref{arbitrary countable graph}, the realisations of which are concrete graphs. (And we could recognise again that from a higher-order perspective, there are many maximally arbitrary countable graphs.) Or, alternatively, we could in a more \emph{in rebus} vein say that graph theory is about all particular graphs.

But there are two reasons why saying that countable graph theory is about the generic countable graph would not fit with mathematical practice. First, this is not how mathematicians talk and think. They do not take the subject of graph theory to be one privileged structure, but instead take graph theory to be about \emph{many} structures. Second, there are graphs that belong to the subject matter of countable graph theory but that are not among the states that the generic countable graph, as described in section \ref{arbitrary countable graph}, can be in. ``\emph{The} unconnected two-element graph'', for instance, is among the graphs that countable graph theory is about \cite{Leitgeb and Ladyman 2007}. But this is not a concrete system: it is not a \emph{particular} graph on $V$.

Nonetheless, there is a natural sense in which it can be said that there are many \emph{arbitrary} countable graphs on an underlying countable collection $B$. By now it is clear how the story goes: \emph{Construct a long repeated list in which each line enumerates all concrete graphs on $B$. Then an arbitrary graph is a thread through this list.}
Among these arbitrary graphs, there are \emph{specific} graphs, such as ``the unconnected two-element graph''. Here by a specific graph I mean a maximally generic system that in one possible state is realised by a particular graph that is isomorphic with every concrete graph that realises it in any other possible situation. (So, in my terminology, \emph{specific} graphs, being generic systems, are not the same as \emph{particular} graphs on $V$, being concrete systems.) 

%Just as one might say that number theory is about the specific numbers, one might say that countable graph theory is about the specific graphs. 

Now it is natural to say ---or so I maintain--- that countable graph theory is about all arbitrary countable graphs. Alternatively, just as one might say that number theory is about the specific numbers, it is natural to say that countable graph theory is about  all specific countable graphs. (Notice that both arbitrary and specific graphs are abstract entities.) 
%Again I will not here try to decide which of these two views is preferable.

There is a higher-order arbitrary object, the realisations of which are arbitrary or specific graphs (the latter of which are themselves also arbitrary systems). In some diluted sense one might say that countable graph theory is about this higher-order arbitrary object. But we should not do so. The reason is that, like the generic countable graph of section \ref{arbitrary countable graph}, the higher-order arbitrary object is not really a \emph{structure}.

Nonetheless, there are natural structural relations between the higher-order object and its realisations. If one were to incorporate these structural relations in the notion of the higher-order object, then one would obtain a notion that is close to what is called a \emph{Lawvere theory} \cite{Lawvere 1964}, \cite{Pettigrew 2011}. Indeed, in his doctoral dissertation Lawvere set out to specify a category that is a generic group (or a generic ring or a generic field). And that category would be defined in such a way that then any specific group (or ring or field) can be identified with a functor from that generic group (or ring or field) into the category of sets.\footnote{I cannot pursue this relation with Category Theory further here, but I am grateful to Richard Pettigrew for drawing my attention to the connection.}

\subsection{Generic structures}\label{gen-struc2}

Mathematical theories are about mathematical structures and their objects. But what mathematical structures and mathematical objects are, is a difficult philosophical question.

On the proposed account, mathematical structures are \emph{generic systems} and mathematical objects are \emph{arbitrary objects}. 
%I call this position \emph{generic structuralism}.

Generic systems are governed by an \emph{instantiation} relation (the relation of being in a state). Therefore my view is a \emph{structuralist} position in the philosophy of mathematics. Clearly it is a Platonistic form of structuralism.

On the proposed view, mathematical structures are more than \emph{mere abstract patterns}, to use an expression from \cite{Resnik 1997}. They differ from Shapiro's ante rem structures. Fine's arbitrary objects play an in-eliminable role in my account, for even \emph{specific} numbers are arbitrary objects and \emph{specific} graphs (or groups or\ldots) are arbitrary systems.

What drives my view, and what has perhaps not been sufficiently recognised by Shapiro's places-objects view, is that in the abstractive movement from system to structure, the notion of object is transformed along with that of system. The appropriate slogan is: \emph{from concrete system to generic system, and from concrete object to arbitrary object}. 

I have tried to show how set theoretical models can give us a feeling of what a generic system such as the generic $\omega$-sequence is like. At the same time, I insisted that the generic $\omega$-sequence is a \emph{sui generis} entity that cannot ontologically be reduced to any pure or impure set.

The previous section contains a sketch of an answer to the question  where the objects out of which the underlying \emph{particular} $\omega$-sequences are constructed come from. There is no plausible way, as far as I can see, that the elements of the underlying plurality of objects can be conjured up, in neo-logicist vein, out of logical air. They must belong to a realm that has ontological priority over the generic  $\omega$-sequence. This realm \emph{can} be taken to be the iterative hierarchy  of pure sets. This hierarchy may be thought of as \emph{reduced} to the generic hierarchy of section \ref{sets}, but it can also be thought of as an alternative that can serve the same purpose as the set theoretic hierarchy.\footnote{After all, the track record of ontological reduction in philosophy is not great.}

However we think of this, the models in sections \ref{arbitrary omega sequences} and \ref{non-alg} are over-simplified: no single countable infinity of sets should be privileged. But this can easily be taken on board: we just have to make our list $\sf{L}$ (much) longer. Take for instance the generic $\omega$-sequence as modelled in section \ref{arbitrary omega sequences}. Let the entries in our revised list $\sf{L}$ consist of exactly \emph{all} the $\omega$-sequence orderings of \emph{all} the countably infinite pluralities of pure sets (or elements of the generic hierarchy), and then proceed as before. A similar story can of course be told for other generic systems.

The fact remains that the ambient domain is given special treatment: it is not itself a \emph{generic} system. It seems to me that this exception can be well motivated along the lines of \cite[p.~144]{Burgess 2015}.\footnote{I should mention that Burgess advocates a form of eliminative rather than non-eliminative structuralism.} Briefly, the motivations from mathematical practice for adopting a structuralist position for particular branches of mathematics (number theory, group theory, analysis, topology,\ldots) just do not seem to extend to foundational theories.

%%%%%%

\section{Rival accounts}\label{comp}

%The Aristotelian reduction of mathematical structures to sets or systems of non-mathematical objects, together with its wholesale elimination of mathematical objects, leaves me unsatisfied. Shapiro's Platonist structuralism fares better, but it contains at least an unclarity at a crucial juncture, i.e., in its stance on the existence and nature of mathematical objects. It is on this point that I have tried to make some philosophical progress.

I now turn to a comparison of my view of mathematical structure with rival accounts of structure and arbitrariness in mathematics. The aim is not so much to show that my account gives better answers to searching philosophical questions as to clarify my stance on some of the issues.

%%%%%%%%%%

\subsection{The incompleteness of mathematical objects}

It is well-documented that  both eliminative and non-eliminative structuralism carry considerable ontological and / or modal commitments. This is particularly so for versions of non-eliminative structuralism. The account that I propose must simply accept such commitments: it is ontologically committed to the existence of arbitrary objects and generic systems. So there is no question of avoiding platonistic commitments. Rather, the question is whether the proposed account gives plausible answers to vexing philosophical problems.

Algebraic theories are about classes of structures that themselves form a more general structure; non-algebraic theories are about particular structures. Nonetheless, particular structures are in a sense more fundamental \cite[p.~1--2]{Isaacson 2011}:
\begin{quote}
Mathematical structures are, roughly, of two kinds, particular (e.g.
the natural numbers) and general (e.g. groups). Mathematics for its
first several thousand years was concerned only with particular structures.
Modern mathematics is much more about general structures, but despite this shift, the reality of mathematics turns ultimately on
the reality of particular structures.
\end{quote}

In contrast to Isaacson's structuralism,\footnote{I do not have space here to go into the details of Isaacson's view.} I take non-algebraic mathematical theories to be not just about structures but also about mathematical \emph{objects}, and thus aim to be more faithful to how mathematicians unreflectively tend to view such structures. Moreover, my aim is to give a detailed  account of the \emph{incompleteness} of mathematical objects, and I do so in terms of the ability of mathematical objects to take on particular values. Indeed, the description of the way in which mathematical objects are incomplete forms the heart of the proposed theory of mathematical structure. 

Shapiro also takes structures to contain objects since an ante rem structure can be taken to be a system. But this move results in mathematical objects not having the right kind of incompleteness, as Hellman's permutation argument shows.\footnote{See section \ref{mathstruc} above.}

The theory of mathematical structure that is advocated in this article does not succumb to the permutation argument. It is clear that the exact analogue of the problem that Shapiro faces does not pose itself. A generic structure is not a state that that same generic structure can be in. Applied to arithmetic, for instance, this means that the generic $\omega$-sequence is not itself an $\omega$-sequence. Nonetheless, one may ask, whether this gets to the heart of the matter.\footnote{Thanks to James Studd for pressing the objection.} Cannot the specific number $0$ in the generic $\omega$-sequence `play the role' of the specific number $1$ and \emph{vice versa}, for instance? But the answer to that question is plainly no: in every state, the value taken by the specific number $0$ plays the $0$-role, and the value taken by the number $1$ plays the $1$-role.
This is the main reason why I take my theory of mathematical structure to be preferable over Shapiro's non-eliminative structuralism.

%%%%%

\subsection{Theories and structures}

Nodelman and Zalta have also proposed a form of ante rem structuralism that is not subject to the permutation objection \cite{Nodelman and Zalta 2014}. I will now briefly discuss their position.

Central to Nodelman and Zalta's account is a distinction between two kinds of predication: \emph{exemplification} and \emph{codification}. Exemplification is the form of predication that we are most familiar with. For instance, we have a case of an ordinary object exemplifying a property when we say that Vladimir Putin is the president of Russia. But abstract objects can codify predicates or groups of predicates. For instance, the abstract object redness codifies the property of being red.

On Nodelman and Zalta's account, mathematical structures are obtained from mathematical theories, along roughly the following lines.
There is a canonical way of constructing properties out of propositions via $\lambda$-abstraction. A mathematical theory can be identified with the collection of mathematical propositions that it logically entails. Therefore, a mathematical theory is associated with the collection of properties constructed out of the propositions that it entails. And there will be a unique abstract object that encodes exactly those properties. So we may \emph{identify} a mathematical structure with this abstract entity \cite[p.~51]{Nodelman and Zalta 2014}.

This is an ante rem form of structuralism that differs from the accounts that we have so far discussed. It associates a unique ante rem structure even with each algebraic theory, and thus comes closer to the account that I propose than to Shapiro's version of non-eliminative structuralism. Moreover, it aims at attributing exactly the right kind of incompleteness to mathematical structures. Consider for instance the theory of dense linear orderings. It will be, on Nodelman and Shapiro's account, about an abstract entity (structure) that encodes all the properties of this theory. In particular, it encodes \emph{neither} countability or uncountability, since there are both countable and uncountable dense linear orderings.

Mathematical structures contain mathematical objects, on Nodelman and Zalta's account. The objects that a structure contains are extracted from the theory from which the structure is obtained. Roughly, objects correspond to the \emph{terms} that the theory contains. The object corresponding to a term of a theory will again be an abstract object that codifies the properties that the theory attributes to it. Thus the Peano Arithmetic structure will be, in some sense, about all and only the familiar objects $0,1,2,\ldots$ And this will have as a consequence that there will be no cross-theory identification of mathematical objects: the number 0 of Peano Arithmetic  will be distinct from the number 0 of Real Analysis.\footnote{So on this point Nodelman and Zalta's account is closer to \cite{Resnik 1997} than to \cite{Shapiro 1997}.}

However, since mathematical objects are identified by means of clusters of properties satisfied by the denotation of a term according to a given theory, there is no room in their account for numerically distinct but strongly indistinguishable objects \cite[p.~73]{Nodelman and Zalta 2014}:
\begin{quote}
An \emph{element} of a structure must be uniquely characterizable in terms of the relations of the structure.
\end{quote} 
But this is in a tension with ways in which mathematicians tend to speak. Take for instance the theory of countable dense linear orderings without endpoints. Intuitively, one would say that any corresponding structure contains countably many numerically distinct objects that are all indistinguishable from each other. But on Nodelman and Zalta's theory, the structure that this theory describes does not contain any objects. 

The moral of this is, I think, that there is more to the objects of a structure than the discriminating powers of our best theory about the structure. On my account, there is a sense in which \emph{even arithmetic} is about many mathematical objects (arbitrary numbers) that are mutually highly indistinguishable from each other.

%%%%%%%%%%

\subsection{Reference and dependence}\label{magidor}

Recall our mathematician from section \ref{arbitrary objects} who enters the lecture room and says ``Let $m$ be an arbitrary natural number''. What is she referring to?

On Breckenridge \& Magidor's account, she is referring to a \emph{specific} natural number but it is in principle impossible for anyone to know which one \cite{Breckenridge and Magidor 2012}. So any arbitrariness lies on the side of the reference relation rather than on the side of the objects referred to. 
Fine in contrast holds that the mathematician is referring to \emph{the (unique) independent arbitrary natural number} \cite[p.~18]{Fine 1985}. Later on in the lecture, Fine would add, she may introduce other arbitrary natural numbers that \emph{depend on} $m$. 

Like Fine's account, the view developed in this article imputes arbitrariness to the objects referred to. But it is nonetheless very different from Fine's. On my account, there are no \emph{independent} arbitrary numbers. All we can say is that all arbitrary numbers are correlated with each other in complicated ways; there is no ontological priority of some over others.

If Fine is right, there is only one natural number that the mathematician's use of `$m$' can refer to: there is then no mystery about how the reference relation finds its target. Nonetheless, as Breckenridge \& Magidor observe, the mathematician could go on to say "Now let $n$ also be an arbitrary natural number." It is not clear in which sense $n$ depends on $m$ (or vice versa)  \cite[p.~391]{Breckenridge and Magidor 2012}.

I conclude from this with Breckenridge \& Magidor that if we take, along with Fine, arbitrary numbers seriously, then there is more than one candidate arbitrary natural number for our mathematician to be referring to when she says "Let $m$ be a natural number". Indeed, we have seen that according to my account, there are \emph{many} completely arbitrary natural numbers. But then the question arises: in virtue of what does `$m$' refer to one arbitrary natural number rather than to another one?

When Breckenridge \& Magidor claim that `$m$' refers to some \emph{specific} natural number, they are faced with a similar question. In response, they say that the reference of `$m$' is a brute, unexplained semantic fact. There is no explanation of why the mathematician refers to the number $172$ (if she does), rather than to the number $3$, for instance. This forces them to deny that reference supervenes on language use \cite[p.~380]{Breckenridge and Magidor 2012}. 

The generic structuralist can take a similar line, and hold that the reference of `$m$' to some \emph{arbitrary} number is a brute semantic fact. But denying the supervenience of reference on language use strikes many as a radical claim.\footnote{This claim is defended in \cite{Kearns and Magidor 2012}.} It might be more reasonable to hold that just as it is \emph{indeterminate} whether the term `$i$' refers to one place in the structure of the complex plane rather than to another,\footnote{The same holds for the term `$-i$'; but it \emph{is} determinate that `$i$' and `$-i$' do not refer to the same place in the structure.} it is indeterminate to which completely arbitrary number the term `$m$' refers.

%%%%%%%%%

\section{In closing}

I have explored a new account of mathematical structure. All I can hope to have achieved is to have articulated and defended this view to the extent that readers are convinced that it deserves a hearing in the community of philosophers of mathematics.

The focus of this article has been mainly on matters metaphysical, although some attention was also payed to aspects of semantic reference. Implications for mathematical epistemology have not been touched upon at all. In particular, Benacerraf's notorious \emph{access problem} \cite{Benacerraf 1970}, which can in some shape or form be raised for every form of mathematical structuralism, was not addressed. This is not because I do not regard it as an important problem; I just do not have space to discuss it in this article.

There are also many technical questions that arise naturally and need to be addressed, such as: What is a suitable formal framework for reasoning about generic structures? Can we formally define notions of indistinguishability of arbitrary mathematical objects in a structure? These questions and others like it are addressed in \cite{Horsten and Speranski 2018}.

%Even though generic structuralism takes non-algebraic theories to be about objects, this account is not vulnerable to the permutation problem. I have observed in section \ref{non-alg} that the ante rem structure of the natural numbers can be seen as a \emph{fusion} of an $\omega$-sequence of arbitrary objects. A non-trivial permutation can be applied to this sequence of arbitrary objects, and this gives rise to a different fusion. But this does not affect the ante rem structure that generic structuralism takes arithmetic to be about, for this new fusion is merely a re-labelling of the generic system that we started out with.

%To conclude, I hope that the perspective from the theory of arbitrary object on the questions of mathematical structuralism will open \emph{new} philosophical avenues, such as bringing probabilistic concepts to bear on discussions about mathematical structuralism.

%%%%%%%%%%%%%%%

\end{document}